\documentclass[10pt,twoside]{article}
\usepackage{graphicx}
\usepackage{amsmath}
\usepackage{Latex-document}

\title{\bf A Numeraire-free and Original \vskip -2mm Probability
 Based Framework for \vskip -2mm Financial Markets\thanks{The work was supported by the
973 project on mathematics of the Ministry of Science and
Technology and the knowledge innovation program of the CAS. The
author wishes to thank Dr. Jianming Xia for helpful
comments.}\vskip 6mm}

\author{Jia-An Yan\vspace*{-0.5cm}\thanks{Academy of Mathematics \& System Sciences,
Chinese Academy of Sciences, Beijing 100080, China. E-mail:
jayan@mail.amt.ac.cn}}
\date{\vspace{-8mm}}

\begin{document}
\maketitle

\thispagestyle{first} \setcounter{page}{861}

\begin{abstract}\vskip 3mm
In this paper, we introduce a numeraire-free and original
probability based framework for financial markets. We reformulate
or characterize fair markets, the optional decomposition theorem,
superhedging, attainable claims and complete markets in terms of
martingale deflators, present a recent result of Kramkov and
Schachermayer (1999, 2001) on portfolio optimization and give a
review of utility-based approach to contingent claim pricing in
incomplete markets.

\vskip 4.5mm



\end{abstract}

\vskip 12mm

\section{Introduction} \label{section 1}\setzero

\vskip-5mm \hspace{5mm}

A widely adopted setting for ``arbitrage-free" financial markets is as follows: one models the price dynamics of
primitive assets by a vector semimartingale, takes the saving account (or bond) as numeraire, and assumes that
there exists an equivalent local martingale measure for the deflated price process of assets. According to the
fundamental theorem of asset pricing (FTAP, for short), due to Kreps (1981) and Delbaen and Schachermayer (1994)
if the deflated price process is locally bounded, this assumption is equivalent to the condition of ``no free
lunch with vanishing risk" (NFLVR for short). However, the property of NFLVR is not invariant under a change of
numeraire. Moreover, under this setting, the market is ``arbitrage-free" only for admissible strategies, the
market may allow arbitrage for static trading strategies with short-selling, and a pricing system using an
equivalent local martingale measure may not be consistent with the original prices of some primitive assets. In
order to remedy these drawbacks, Yan (1998) introduced the numeraire-free notions of ``allowable strategy" and
fair market.  In this paper, we will further present a numeraire-free and original probability based framework for
financial markets in a systematic way.

The paper is organized as follows: In Section 2, we introduce the
semimartingale model, define the notion of martingale deflator. In
Sction 3, we reformulate Kramkov's optional decomposition theorem
in terms of martingale deflators, and give its applications to the
superhedging of contingent claims and the characterizations of
attainable claims and complete markets. In Section 4, we present a
recent result of Kramkov and Schachermayer (1999, 2001) on optimal
investment and give a review of utility-based approach to
contingent claim pricing in incomplete markets.

\section{Semimartingale model and basic concepts}\label{section 2}\setzero
\vskip-5mm \hspace{5mm}

We consider a security market model in which the uncertainty and information structure are described by a
stochastic basis $(\Omega,{\cal F},P; ({\cal F}_t))$ satisfying the usual conditions with ${\cal F}_0$ being
trivial. We call $P$ the original (or objective) probability. It models the ``real world" probability.

The market consists of $d$ (primitive) assets whose price
processes $(S^i_t), i=1,\cdots,d$ are assumed to be non-negative
semimartingales with initial values non-zero. We further assume
that the process $\sum_{i=1}^dS^i_t$ is strictly positive and that
each $S^i_t$ vaniches on $[T^i,\infty)$, where
$T^i(\omega)=\inf\{t>0: S^i_t(\omega)=0, \ {\rm or}\ \
S^i_{t-}(\omega)=0\}$ stands for  the ruin time of the company
issuing asset $i$. We will see later that this latter assumption
is automatically satisfied for a fair market, since any
non-negative supermartingale satisfies this property. In the
literature, it was assumed that all primitive assets have strictly
positive prices.

Let $S_t=(S^1_t,\cdots, S^d_t)$. Throughout the paper, we will use
the following notation:
$$S^*_t=\left(\sum_{i=1}^dS^i_0\right)^{-1}\sum_{i=1}^dS^i_t.$$
By assumption, $S^*_t$ is a strictly positive semimartingale. In
the literature on mathematical finance, one often takes a
primitive asset whose price never vanishes as  numeraire. In our
model, such a primitive numeraire asset may not exist. However, by
our assumption on the model,  we can always take $S^*_t$ as
numeraire.

\subsection{Self-financing strategy}

\vskip-5mm \hspace{5mm}

A {\it trading strategy} is an $R^d$-valued ${\cal F}_t$-predictable process $\theta(t)=(\theta^1(t),\cdots$,
$\theta^d(t))$, which is integrable w.r.t. the semimartingale $S_t$. Here $\theta^i(t)$ represents the numbers of
units of asset $i$ held at time $t$. The wealth $W_t(\theta)$ at time $t$ of a trading strategy $\theta$ is
$W_t(\theta)=\theta(t)\cdot S_t,$ where $a\cdot b$ denotes the inner product of two vectors $a$ and $b$. A trading
strategy $\theta$ is said to be {\it self-financing}, if
\begin{equation}\label{2.1}
 W_t(\theta)=W_0(\theta)+\int_0^t\theta(u)dS_u.
\end{equation}
In this paper we use notation $\int_0^tH_udX_u$ or $(H.X)_t$ to
denote the integral of $H$ w.r.t. $X$  over the interval $(0,t]$.
In particular, we have $(H.X)_0=0$.

The following theorem concerns a result on stochastic integrals of
semimartingales, which represents an important property of
self-financing strategies. It was given in Xia and Yan (2002).

{\bf Theorem 2.1} \it Let $X$ be an $R^d$-valued semimartingale
and $H$ an $R^d$-valued predictable process. If $H$ is integrable
w.r.t. $X$ and
\begin{equation}\label{2.2}
H_t\cdot X_t=H_0\cdot X_0+\int_0^tH_sdX_s,
\end{equation}
 then for
any real-valued semimartingale $y$, $H$ is integrable w.r.t. $yX$
and
\begin{equation}\label{2.3}
y_t(H\cdot X)_t=y_0 (H\cdot X)_0+\int_0^tH_sd(yX)_s.
\end{equation}\rm

As a consequence of Theorem 2.1, we obtain the following

{\bf Theorem 2.2}\ \it 1) For any given $R^d$-valued
$S$-integrable predictable process $\theta(t)$ and a real number
$x$ there exists a real-valued predictable process $\theta^*(t)$
such that $\{\theta^*(t)1_d+\theta(t)\}$ is a self-financing
strategy with initial wealth $x$, where $1_d$ is the
$d$-dimensional vector $(1,1,\cdots,1)$.

2) A strategy $\theta$ is self-financing if and only if
$d\widetilde W_t(\theta)=\theta(t) d\widetilde S_t$, where
$\widetilde S_t=S_t(S^*_t)^{-1}, \widetilde
W_t(\theta)=W_t(\theta)(S^*_t)^{-1}$. \rm

\subsection{Fair market and fundamental theorem of asset pricing}

\vskip-5mm \hspace{5mm}

Now we consider a finite time horizon $T$. In Yan (1998), we
introduced the notions of allowable strategy and fair market under
assumption that all price processes of assets are strictly
positive. The following definitions extend these notions  to the
present model.

{\bf Definition 2.1} \it A strategy $\theta$ is said to be
allowable, if it is self-financing and there exists a positive
constant $c$ such that the wealth $W_t(\theta)$ at any time $t$ is
bounded from below by $-cS^*_t$.\rm

{\bf Definition 2.2} \it A market is said to be fair if there
exists a probability measure $Q$ equivalent to the original
probability measure $P$ such that the deflated price process
$(\widetilde S_t)$ is a (vector-valued) $Q$-martingale.  \rm

We  call such a $Q$ an {\it equivalent martingale measure} for the
market. Throughout the sequel we denote by ${\cal Q}$ the set of
all equivalent martingale measures.

If the market is fair, the deflated  wealth process of any
allowable strategy is a local $Q$-martingale, and consequently, is
also a $Q$-supermartingale,  for all $Q\in{\cal Q}$.

By the main theorem in Delbaen and Schachermayer (1994), Yan
(1998) obtained an intrinsic  characterization of fair markets.
This result can be regarded as a numeraire-free version of the
FTAP due to Kreps (1981) and Delbaen and Schachermayer (1994). The
same result is valid for our more general model.

{\bf Theorem 2.3}\  \it The market is fair if and only if there is
no sequence $(\theta_n)$ of allowable strategies with initial
wealth 0 such that $W_T(\theta_n)\ge -\frac1n S^*_T$\,a.s., $\
\forall n\ge 1,$, and such that $W_T(\theta_n)$ a.s. tends to a
non-negative random variable $\xi$ satisfying $P(\xi>0)>0$.\rm

{\bf Remark}\ \ If we take $S^*_t$ as numeraire and consider the
market in deflated terms, the condition in Theorem 2.3 is the
NFLVR condition introduced in Delbaen and Schachermayer (1994).

\subsection{Martingale deflators}

\vskip-5mm \hspace{5mm}

In principle, we can take any strictly positive semimartingale as
a numeraire, and its reciprocal as a deflator.

{\bf Definition 2.3} \it A strictly positive semimartingale $M_t$
with $M_0=1$ is called a martingale deflator for the market, if
the deflated price processes $(S^i_tM_t), i=1,\cdots,d$ are
martingales under the original probability measure $P$.\rm

In the literature, such a deflator $M$ is called ``state price
deflator". Here we propose to name it as ``martingale deflator". A
martingale deflator $M$ is uniquely determined by its terminal
value $M_T$. In fact, we have $M_t=(S^*_t)^{-1}E[M_TS^*_T|{\cal
F}_t]$.

In terms of martingale deflators, a market is  fair if and only if
there exists a martingale deflator for the market.

Assume that the market is fair. We denote by ${\cal M}$ the set of
all martingale deflators, and denote by ${\cal Q}$ the set of all
equivalent martingale measures, when $S^*_t$ is taken as
numeraire. Note that there exists a one-to-one correspondence
between  ${\cal M}$ and ${\cal Q}$. If $M\in {\cal M}$, then
$\frac{dQ}{dP}=M_TS^*_T$ define an element $Q$ of ${\cal Q}$. If
$Q\in {\cal Q}$, then we can define an element $M$ of ${\cal M}$
with  $M_T=\frac{dQ}{dP}(S^*_T)^{-1}$. If ${\cal M}$ (or ${\cal
Q})$ contains only one element, the market is said to be complete.
Otherwise, the market is said to be incomplete.

We will see in following sections that the use of  martingale
deflators instead of equivalent martingale measures has some
advantages in handling financial problems.

\section{Optional decomposition theorem and its applications} \label{section 3}\setzero
\vskip-5mm \hspace{5mm}

The optional decomposition theorem of Kramkov is a very useful tool in mathematical finance. It generalizes the
classical Doob-Meyer decomposition theorem for supermartingales. This kind of decomposition was first proved by El
Karoui and Quenez (1995), in which the process involved is the value process of a superheding strategy for a
contingent claim in an incomplete market modelled by a  diffusion process. Kramkov (1996) extended this result to
the general semimartingale setting, but under the assumption that the underlying semimartingale  is locally
bounded and the supermartingale to be decomposed is non-negative and locally bounded. F\"ollmer and Kabanov (1998)
removed any boundedness assumption. But in both papers, the theorem was formulated in the setting that there
exists equivalent local martingale measures for the underlying semimartingales.

\subsection{Optional decomposition theorem in terms of martingale deflators}

\vskip-5mm \hspace{5mm}

Based on Theorem 2.1, Xia and Yan (2002) obtained  the following
version of the optional decomposition theorem in the equivalent
martingale measure setting.

{\bf Theorem 3.1} \it Let $Y$ be a vector-valued semimartingale
with non-negative components. Assume that the set ${\cal Q}$ of
equivalent martingale measures for $Y$ is nonempty. If $X$ is a
local ${\cal Q}$-supermartingale, i.e. local $Q$-supermartingale
for all \, $Q\in{\cal Q}$, then there exist an adapted, right
continuous and increasing process $C$ with $C_0=0$, and a
$Y$-integrable predictable process $\varphi$ such that
$$X=X_0+\varphi.Y-C.$$
Moreover, if $X$ is non-negative, then $\varphi.Y$ is a local
${\cal Q}$-martingale. \rm

The following theorem is a reformulation of Theorem 3.1 in terms of martingale deflators.

{\bf Theorem 3.2} \it Assume that the market is fair. We denote by
${\cal M}$ the set of all martingale deflators. Let $X$ be a
semimartingale. If $XM$ is a local supermartingale for all \,
$M\in{\cal M}$, then there exist an adapted, right continuous and
increasing process $C$ with $C_0=0$, and an $S$-integrable
predictable process $\varphi$ such that
$$X=X_0+\varphi.S-C.$$
Moreover, if $X$ is non-negative, then $(\varphi.S)M$ is a local
martingale for all \, $M\in{\cal M}$. \rm

{\bf Proof} Let $\widetilde S_t=S_t(S^*_t)^{-1}$ and $\widetilde
X_t= X_t(S^*_t)^{-1}$. Let ${\cal Q}$ denote the set of all
martingale measures for $\widetilde S$. Then $\widetilde X$ is a
local ${\cal Q}$-supermartingale. By  Theorem 3.1 we have
$$\widetilde X=X_0+\psi.\widetilde S-D,$$
where $D$ is an adapted, right continuous and increasing process
with $D_0=0$. By Theorem 2.2 there exists a real-valued
predictable process $\theta^*(t)$ such that
$\{\theta^*(t)1_d+\psi(t)\}$ is a self-financing strategy with
initial wealth $X_0$. Since $\sum_{i=1}^d \widetilde
S^i_t=\sum_{i=1}^d S^i_0$, we have $\theta^*(t)1_d.\widetilde
S=0$. Consequently,
$$X_0+((\theta^*1_d+\psi).S)_t=S^*_t(X_0+((\theta^*1_d+\psi).\widetilde S)_t)
=S^*_t(X_0+(\psi.\widetilde S)_t)=X_t+S^*_tD_t.$$ Put
$\varphi=(\theta^*-D_{-})1_d+\psi$ and $C=S^*.D$, we get the
desired decomposition.

\subsection{Superhedging}

\vskip-5mm \hspace{5mm}

By a {\it contingent claim} (or {\it derivative}) we mean a
non-negative ${\cal F}_T$-measurable random variable. Let $\xi$ be
a contingent claim. In general, one cannot find a self-financing
strategy to perfectly replicate $\xi$. It is natural to raise the
question: Does there exist an admissible strategy with the minimal
initial value, called superhedging strategy, such that its
terminal wealth is no smaller than the claim $\xi$? Here and
henceforth, by an {\it admissible strategy} we mean a
self-financing strategy with non-negative wealth process. For a
market with diffusion model, this problem has been solved by El
Karoui and Quenez (1995). For a general semimartingale model, it
was solved by Kramkov (1996) and F\"ollmer and Kabanov(1998) using
the optional decomposition theorem. The initial value of the
superhedging strategy is called the {\it cost of superhedging
$\xi$}. It can be considered as the ``selling price" or ``ask
price " of $\xi$.

In a fair market setting, based on the corresponding result of
Kramkov (1996), Xia and Yan (2002) proved the following result: if
$\sup_{Q\in{\cal Q}}E_Q\left[(S^*_T)^{-1}\xi\right]<\infty$, then
the cost at time $t$ of superhedging the claim $\xi$ is given by
\begin{eqnarray}\label{(3.1)}
U_t=\mbox{esssup}_{Q\in{\cal
Q}}S^*_tE_Q\left[\left.{(S^*_T)^{-1}\xi}\right|{\cal
F}_t\right].\end{eqnarray} $U$ is the smallest non-negative ${\cal
Q}$-supermartingale with $U_T\ge \xi$. In terms of martingale
deflators, we can rewrite (3.1) as
\begin{eqnarray}\label{(3.2)}
U_t=\mbox{esssup}_{M\in{\cal
M}}M_t^{-1}E\left[\left. M_T\xi\right|{\cal
F}_t\right].
\end{eqnarray}

Using the optional decomposition theorem F\"ollmer \& Leukert
(2000) showed that the optional decomposition of a suitably
modified claim gives a more realistic hedging (called {\it
efficient hedging}) of a contingent claim. This result can be also
reformulated in terms of martingale deflators.

\subsection{Attainable claims and completeness of the market}

\vskip-5mm \hspace{5mm}

Xia and Yan (2002) introduced the notions of regular and strongly
regular strategies. We reformulate them in terms of martingale
deflators.

{\bf Definition 3.1} \it  A self-financing strategy $\psi$ is said
to be regular (resp. strongly regular), if for some (resp. for all
) $M\in {\cal M}$,  $W_t(\psi)M_t$ is a martingale. A contingent
claim is said to be  attainable if it can be replicated by a
regular strategy.\rm

By Theorem 3.2, one can easily deduce the following
characterizations  for attainable claims and complete markets.

{\bf Theorem 3.3} \it Let  $\xi$ be a contingent claim  such that
$\sup_{M\in {\cal M}} E \left[{\xi M_T}\right]<\infty.$ Then $\xi$
is attainable (resp. replicatable by a strongly regular strategy)
if and only if the above supremum is attained by  an $M^*\in {\cal
M}$ (resp. $E[M_T\xi]$ doesn't depend on $M\in{\cal M}$).\rm

{\bf Theorem 3.4} \it The market is complete if only if any
contingent claim $\xi$ dominated by $S^*_T$ is attainable, or
equivalently, $E[M_T\xi]$ doesn't depend on $M\in{\cal M}$.\rm

\section{Portfolio optimization and contingent claim pricing} \label{section 4}\setzero
\vskip-5mm \hspace{5mm}

The portfolio optimization and contingent claim pricing and hedging are three major problems in mathematical
finance. In a market where assets prices follow an exponential L\'evy process, the portfolio optimization problem
was studied in Kallsen (2000). In the general semimartingale model, for utility functions $U$ with effective
domains ${\cal D}(U)=R_+$, the portfolio optimization problem  was completely solved by Kramkov and
Schachermayer(1999, 2001), henceforth K-S(1999, 2001). Bellini \& Frittelli(2002) and Schachermayer (2002) studied
the problem for utility functions $U$ with ${\cal D}(U)=R$. The relationship between portfolio optimization and
contingent claim pricing was studied in Frittelli(2000) and Goll \& R\"uschendorf(2001), among others. In what
concerning the problem of hedging contingent claims, we refer the reader to Schweizer (2001) for quadratic
hedging, F\"ollmer \& Leukert (2000) for efficient hedging, and Delbaen et al. (2001) for exponential hedging.

In this section, under our framework, we will present the main
results of K-S(1999, 2001) and give a review of utility-based
approach to contingent claim pricing.

\subsection{Expected utility maximization}

\vskip-5mm \hspace{5mm}

We consider an agent whose objective is to choose a trading
strategy to maximize the expected utility from terminal wealth at
time $T$. In the sequel, we only consider such a utility function
$U:(0,\infty)\longrightarrow R$, which is strictly increasing,
strictly concave, continuously differentiable and satisfies
$\lim_{x\downarrow 0}U^\prime(x)=\infty,
\lim_{x\to\infty}U^\prime(x)=0.$ We denote by $I$ the inverse
function of $U^\prime$. The conjugate function $V$ of $U$ is
defined as
$$V(y)=\sup_{x>0}[U(x)-xy]=U(I(y))-yI(y), \ y>0.$$

For $x>0$, we denote by ${\cal A}(x)$ the set of all admissible
strategies $\theta$ with initial wealth $x$. For $x>0, y>0, $ we
put
$${\cal X}(x)=\{W(\theta): \theta\in {\cal A}(x)\}, \qquad {\cal X}={\cal X}(1),$$
$${\cal Y}=\{Y\ge 0: Y_0=1,\ YX \ {\rm is \ a \ supermartingale} \ \forall
X\in {\cal X}\}, \qquad {\cal Y}(y)=y{\cal Y},$$
$${\cal C}(x)=\{g\in L^0(\Omega, {\cal F}_T, P), 0\le g\le X_T,\ {\rm for \
some \ } X\in{\cal X}(x)\}, \qquad {\cal C}\hat={\cal C}(1),$$
$${\cal D}(y)=\{h\in L^0(\Omega, {\cal F}_T, P), 0\le h\le Y_T,\ {\rm for \
some \ } Y\in{\cal Y}(y)\}, \qquad {\cal D}\hat={\cal D}(1).$$ The
agent's optimization problem is:
$$\widehat \psi(x)={\rm arg}\max_{\psi\in{\cal A}(x)}E\left[U(W_T(\psi))\right].$$
To solve this problem we consider two optimization problems (I)
and (II):
$$\widehat X(x)={\rm arg}\max_{X\in{\cal
X}(x)}E\left[U(X_T)\right];\ \ \ \ \widehat Y(y)={\rm
arg}\min_{Y\in{\cal Y}(y)}E\left[V(Y_T)\right].$$ Problem (II) is
the dual of problem (I). Their value functions  are
$$u(x)=\sup_{X\in{\cal
X}(x)}E\left[U(X_T)\right], \qquad v(y)= \inf_{Y\in{\cal
Y}(y)}E\left[V(Y_T)\right].$$

The following theorem is the reformulation of the main results of
K-S(1999, 2001) under our framework.

{\bf Theorem 4.1} \it Assume that there is a $\psi\in {\cal A}(1)$
such that $W_T(\psi)\ge K$ for a positive constant $K$ (e.g.,
$S^*_T\ge K$). If $v(y)<\infty, \forall y>0$, then the value
functions $u(x)$ and $v(y)$ are conjugate in the sense that
$$v(y)=\sup_{x>0}[u(x)-xy], \ \ \ u(x)=\inf_{y>0}[v(y)+xy], $$
and we have:

1. For any $x>0$ and $y>0$,  both optimization problems (I) and
(II) have unique solutions $\widehat X(x)$ and $\widehat Y(y)$,
respectively.

2. If $y=u'(x)$, then $\widehat X_T(x)=I(\widehat Y_T(y))$ and the
process $\widehat X(x)\widehat Y(y)$ is a martingale.

3. $v(y)=\inf_{M\in{\cal M}} E\left[V(yM_T)\right]$. \rm

{\bf Proof} \  The proof is almost the same as that in K-S (1999,
2001). We indicate below main differences from K-S (1999, 2001).
Obviously, ${\cal C}$ and ${\cal D}$ are convex sets. By
Proposition 3.1 and a slight modification of Lemma 4.2 in
K-S(1999), one can show that ${\cal C}$ and ${\cal D}$ are closed
under the convergence in probability. For Items 1 and 2, as in
Lemma 3.2 of K-S(1999) and Lemma 1 of K-S(2001), in order to prove
the families $(V^-(h))_{h\in{\cal D}(y)}$ and $(U^+(g))_{g\in{\cal
C}(x)}$ are uniformly integrable, we need to use a fact that
${\cal C}$ contains a positive constant. In our case, we have
indeed  $K\in{\cal C}$, since by assumption  $K\le W_T(\psi)$ for
some $\psi\in {\cal A}(1)$. As for Item 3, according to
Proposition 1 in K-S(2001) we only need to show $\widehat{\cal
D}=\{ M_T:\ M\in{\cal M}\}$ satisfies the following conditions:
\begin{itemize}
  \item For any $g\in{\cal C}$,
  $\sup_{h\in\widehat{\cal D}}E[gh]=\sup_{h\in{\cal D}}E[gh]$
  \item $\widehat{\cal D}\subset{\cal D}$, $\widehat{\cal D}$ is convex and closed under
  countable convex combinations.
  \end{itemize}
The first condition follows easily from (3.2), the second one is
trivial.

\subsection{Utility-based approach to contingent claim pricing}

\vskip-5mm \hspace{5mm}

Assume that the market is fair. Let $\xi$ be a contingent claim
such that $M_T\xi$ is integrable for some $M\in {\cal M}$. We put
\begin{eqnarray}\label{(4.1)}
V_t=(M_t)^{-1}E \left[{M_T}\xi\,\vert\,{\cal
F}_t\right].\end{eqnarray} If we specify $(V_t)$ as the price
process of an asset generated by $\xi$, then the market augmented
with this derivative asset is still fair, because $M$ is still a
martingale deflator for the augmented market.  So we can define
$(V_t)$ as a ``fair price process" of $\xi$. This pricing rule is
consistent with the original price processes of primitive assets.
However, if the market is incomplete (i.e., the martingale
deflator is not unique) we cannot, in general, define uniquely the
fair price process of a contingent claim.

In deflated terms, pricing of contingent claims in an incomplete
market consists in choosing  a reasonable martingale measure.
There are several approaches to make such a choice. A well-known
one is the so-called ``utility-based approach". The basic idea of
this approach is as follows. Assume that the representative agent
in the market has preference represented by a utility function. In
certain cases, the dual optimization problem (II) may produce a so
called minimax martingale measure (MMM for short).

Now under our framework we show how the expected utility
maximization problem is linked by duality to a martingale
deflator. Assume that the solution ${\widehat Y}(y)$ of the dual
optimization problem (II) lies in $y{\cal M}$. We put $\widehat
M(y)=y^{-1}{\widehat Y}(y)$. Then $\widehat M(y)\in {\cal M}$, and
we have
$$\widehat M(y)={\rm arg}\min_{M\in{\cal M}}E\left[V(yM_T)\right].$$
We call $\widehat M(y)$ the {\it minimax martingale deflator}.

The following theorem gives a necessary and sufficient condition
for the existence of the minimax martingale deflator.

{\bf Theorem 4.2} \it Assume that there is a $\psi\in {\cal A}(1)$
such that $W_T(\psi)\ge K$ for a positive constant $K$ (e.g.,
$S^*_T\ge K$), and that $v(y)<\infty$ for all $y>0$. Let $x>0$ be
the agent's initial wealth  and $M^*\in {\cal M}$. In order that
$M^*\in {\cal M}$ is the minimax martingale deflator corresponding
to the utility function $U$ if and only if there exist $y>0$ and
$X^*\in {\cal X}(x)$ such that $X^*_T=I(yM^*_T)$ and
$E[M^*_TX^*_T]=x$. If it is the case, then $X^*$ solves the
optimization problem (I). \rm

{\bf Proof} We only need to prove the sufficiency of the
condition. We have the following inequality
$$U(I(z))\ge U(w)+z[I(z)-w], \ \ \forall w>0, z>0.$$
If we replace $z$ and $w$ by $yM^*_T$ and $X_T\in {\cal X}(x)$ and
take expectation w.r.t. $P$, we get immediately that
$E[U(X^*_T)]\ge E[U(X_T)]$ for all $X\in {\cal X}(x)$. This shows
that $X^*$ solves the optimization problem (I). On the other hand,
since $X^*_T=I(yM^*_T)$ and the assumption $E[M^*_TX^*_T]=x$
implies that $M^*X^*$ is a martingale, by Theorem 3.1, $yM^*$ must
solve the optimization problem (II). In particular, $M^*$ is the
minimax martingale deflator.

Now assume the minimax martingale deflator $\widehat M(y)$ exists.
Let $\xi$ be a contingent claim. If we use $\widehat M(y)$ to
compute a fair price of $\xi$ by (4.1), then it coincides with the
fair price of Davis (1997), which is derived through the so-called
``marginal rate of substitution" argument. In fact, the Davis'
fair price of $\xi$ is defined by
$$\hat\pi(\xi)=\frac{E[U'(\widehat X_T(x))\xi]}{u'(x)}.$$
Since $y=u'(x)$ and $U'(\widehat X_T(x))=\widehat Y(y)$, we have
$\hat\pi(\xi)=E[\widehat M_T(y)\xi]$.

Now we explain the economic meaning of Davis' fair price of a
contingent claim. Let $\xi$ be a contingent claim with $E[\widehat
M_T(y)\xi]<\infty$. Put $\xi_t=(\widehat M_t(y))^{-1}E[\widehat
M_T(y)\xi |{\cal F}_t]$. We augment the market with derivative
asset $\xi$, and consider the portfolio maximization problem in
the new market. Then it is easy to see that $\widehat Y(y)$ is
still the solution of the dual optimization problem (II) in the
new market. Consequently, the value function $v$ and its conjugate
function $u$ remain unchanged. By Theorem 4.1, $\widehat X_T(x)$
solves again the optimization problem (I) in the new market. This
shows that if the price of a contingent claim is defined by Davis'
fair price, no trade on this contingent claim increases the
maximal expected utility in comparison to an optimal trading
strategy. This fact was observed in Goll and R\"uschendorf (2001).

Note that in general the MMM (or minimax martingale deflator)
depends on the agent's initial wealth  $x$. This is a disadvantage
of the utility-based approach to contingent claim pricing.
However, for utility functions  $\ln x, \frac {x^p}p, -e^{- x}$,
where $p\in (-\infty, 1)\setminus \{0\}, \alpha>0$, the MMM is
independent of the agent's initial wealth $x$. This is due to the
fact that the conjugate functions of the above utility functions
are $-\ln x-1, -\frac{p-1}px^{\frac p{p-1}}, -x+x\ln x$,
respectively, and that $E[{dQ}/{dP}]=1$ for any equivalent
martingale measure $Q$. Under our framework, the situation is a
little different: for exponential utility function $U(x)=-e^{-x}$,
the minimax martingale deflator depends still on the agent's
initial wealth $x$.

For $U(x)=-e^{-x}$, the corresponding MMM is called the {\it
minimal entropy martingale measure}. We refer the reader to
Frittelli (2000), Miyahara (2001) and Xia \& Yan (2000) for
studies on the subject. If $U(x)=\ln x$, the minimax martingale
deflator $\widehat M$, if it exists, is nothing but the reciprocal
of the wealth process $\widehat X(1)$ of the growth optimal
portfolio. Yan, Zhang \& Zhang (2000) worked out explicit
expressions for growth optimal portfolios in markets driven by a
jump-diffusion-like process or by a L\'evy process. See also
Becherer (2001) for a study on the subject.

\section{Concluding remarks} \label{section 5}\setzero
\vskip-5mm \hspace{5mm}

We have introduced a numeraire-free and original probability based framework for financial markets. This framework
has the following advantages: Firstly, it permits us to formulate financial concepts and results in a
numeraire-free fashion. Secondly, since the original probability models the ``real world" probability, one can
investigate the martingale deflators by statistical methods using market data. Thirdly, using martingale deflators
to deal with problems of pricing and hedging as well as portfolio optimization is sometimes more convenient than
the use of equivalent martingale measures. Lastly, our framework includes the traditional one with deflated terms
as a particular case. In fact, if the price process of one primitive asset is the constant 1, our framework is
reduced to the traditional one.

\label{lastpage}

\end{document}